\documentclass[12pt]{article}
\usepackage{cite,makeidx,amscd,latexsym,shapepar,hyperref,amsmath,amsthm,amsopn,amstext,amsfonts,amsbsy}
\usepackage{float}
\usepackage{graphicx}

\begin{document}

\author{ Alireza Khalili Golmankhaneh $^a$\\Dumitru Baleanu $^b$
\footnote{Tel:+903122844500, ~~~~Fax:+903122868962~\newline
\textit{E-mail addresses}:~dumitru@cankaya.edu.tr~}
\\$^a$\textit{Department of Physics,~Islamic Azad University, Urmia
Branch,}\\\textit{PO Box 969,~Urmia,~Iran}\\
$^{b}$\textit{Department of Mathematics and Computer Science}\\
\textit{$\c{C}ankaya$} \textit{University, 06530 Ankara,
Turkey}\\\textit{and}\\ \textit{Institute of Space Sciences,}\\
\textit{P.O.BOX, MG-23, R76900,~Magurele-Bucharest, Romania}  }

\title{\textbf{New Derivatives on  Fractal Subset of Real-line }}\maketitle \large
%========================================================

\begin{abstract}
In this manuscript we introduced the generalized fractional
Riemann-Liouville and Caputo like derivative for functions defined
on fractal sets. The Gamma, Mittag-Leffler and Beta functions were
defined on the fractal sets. The non-local Laplace transformation
is given and applied for solving linear and non-linear fractal
equations. The advantage of using these new nonlocal derivatives
on fractals subset of real-line lies in the fact that they are
used for better modelling of processes with memory effect.
\end{abstract}


\section{Introduction}

The calculus involving arbitrary orders  of derivatives and
integrals is called the fractional calculus.~Recently, the
fractional calculus has found many applications in several areas
of science and engineering
\cite{book-1,book-2,book-b3,book-b4,book-b6}. The nonlocal
property of the fractional derivatives and integrals is used to
model the processes  with memory effect \cite{book-1,book-2}. For
example, the fractional derivatives are used to model more
appropriate  the dynamics of the non-conservative systems in the
Hamilton, Lagrange and Nambu mechanics
\cite{book-3,book-4,book-b4-4,book-5,book-6}.~The continuous  but
non-differentiable functions admit the  local fractional
derivatives \cite{book-7}. The local fractional derivative give
the measure on fractal sets. Consequently, recently the
$F^{\alpha}$-calculus  on the fractal subset of real line and
fractal curves is built as a framework \cite{book-8,book-9}.
Fractal analysis is established by many researchers by using
different methods \cite{book-10,book-11,book-12,book-13}. Using
$F^{\alpha}$-calculus the Newton,  Lagrange and Hamilton mechanics
were built on fractal sets   \cite{book-14,book-15}. Also, the
Schr\"{o}dinger's equation on fractal curve was derived in
\cite{book-16,book-17,book-1722}. Motivated by the above mentioned
interesting results,  in this work, we define the non-local
derivative on fractal sets. These new  derivatives can be
successfully  used to derive new mathematical models on fractal
sets involving processes
with memory.\\

We organize our manuscript as follows:\\

In Section \ref{se-2}, we gave a brief exposition of
$F^{\alpha}$-calculus and defined fractal Gamma and Beta
functions. In Section \ref{se-3} we defined the non-local
derivative on fractals as a generalized Riemann-Liouville and
Caputo fractional derivatives.~In Section \ref{se-4},
Mittag-Leffler function and non-local Laplace fractional on
fractal sets are introduced. We solved the non-local differential
equations on fractal using the suggested methods. Section
\ref{se-5} was devoted to our conclusion.
\section{ A review of fractional local derivatives \label{se-2}}

In this section, we review the $F^{\alpha}$-calculus
\cite{book-8,book-9}.
\subsection{Calculus on fractal subset of real-line}

\begin{figure}[H]
\center
   \includegraphics[scale=0.6]{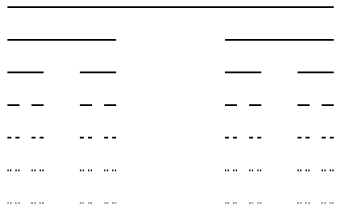}
 \caption{The finite iteration for constructing  the triadic Cantor set.}\label{cantor1f}
\end{figure}

\begin{figure}[H]
\center
   \includegraphics[scale=0.6]{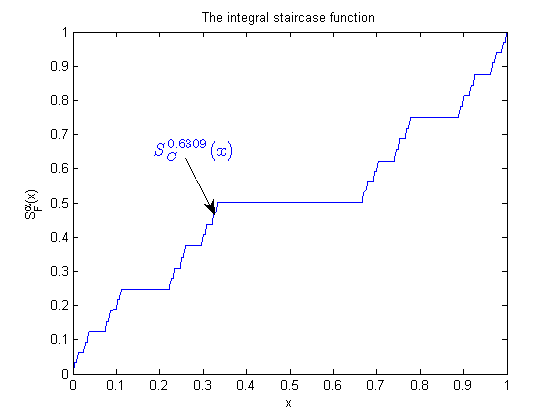}
 \caption{We plot the integral staircase function for triadic Cantor.}\label{hddho12}
\end{figure}

The fractal geometry is the geometry of the real world
\cite{book-1}. Fractal shape is the object with the fractional
dimension and the self similarity property \cite{book-5,book-6}.
In the seminal paper,  Parvate and Gangal have established a
calculus on fractals which similar to Riemann integration. The
suggested framework become a mathematical model for many phenomena
in fractal media \cite{book-8,book-9}. We recall that the triadic
Cantor set is a fractal that can be obtained by an iterative
process. In figure \ref{cantor1f} we show the Triadic Cantor set
\cite{book-10}.
The integral staircase function for the triadic Cantor set is
defined  as \cite{book-8,book-9}.
\begin{equation}
S_{F}^{\alpha}(x)=\left\{
                    \begin{array}{ll}
                      \gamma^{\alpha}(F,a,b), & ~\textmd{if}~ x\geq a_{0},  \\
                      -\gamma^{\alpha}(F,a,b), &~ \textmd{othewise}.
 \end{array}
                  \right.
               \end{equation}
where $\alpha$ is the $\gamma$-dimension of triadic Cantor set. In
Figure \ref{hddho12} we plot the integral staircase function for
triadic Cantor set.
The definitions of  $F^{\alpha}$-limit, $F^{\alpha}$-continuity
and $F^{\alpha}$-integration are given  in the ref.
\cite{book-8,book-9}. The $F^{\alpha}$-differentiation is denoted
by $D_{F}^{\alpha}$ and it is defined as
\begin{equation}\label{jjuh}
  D_{F}^{\alpha}f(x)=\left\{
                       \begin{array}{ll}
                         F-\lim_{y\rightarrow x}\frac{f(y)-f(x)}{S_{F}^{\alpha}
                         (y)-S_{F}^{\alpha}(x)}, &   \textmd{if}~~x \in F,\\
                         0, & \textmd{otherwise,}
                       \end{array}
                     \right.
\end{equation}
if the limit exists \cite{book-8,book-9}.\\
\textbf{Definition 1.} The Gamma function with the fractal support
is defined as
\begin{equation}\label{cderf}
  \Gamma_{F}^{\alpha}(x)=\int_{S_{F}^{\alpha}(0)}^{S_{F}^{\alpha}(\infty)}
e^{-S_{F}^{\alpha}(t)}S_{F}^{\alpha}(t)^{S_{F}^{\alpha}(x)-1}d_{F}^{\alpha}t,
\end{equation}
where
\begin{equation}\label{xswazq}
  e^{-S_{F}^{\alpha}(t)}=F-\lim_{n \rightarrow \infty}(1-\frac{S_{F}^{\alpha}(t)}{n})^n.
\end{equation}
\begin{figure}[H]
\center
   \includegraphics[scale=0.5]{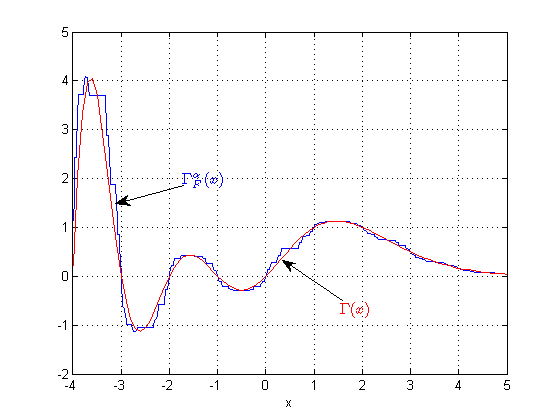}
 \caption{We sketch  the fractal Gamma function which is compared with the standard case.
    }\label{894544}
\end{figure}
\textbf{Definition 2.} The fractal Beta function on the fractal
set is defined as follows
\begin{equation}\label{zaws}
    B_{F}^{\alpha}(r,w)=\int_{0}^{1}S_{F}^{\alpha}(\zeta)^{r-1}
    (1-S_{F}^{\alpha}(\zeta))^{w-1}d_{F}^{\alpha}\zeta,
\end{equation}
which is called two-parameter $(r,w)$ fractal integral, where
$~\Re(r)>0$ and $\Re(w)>0$. \\
In the following we present some  properties of fractal Beta
function.\\
1) The fractal Beta function has a symmetry
$B_{F}^{\alpha}(w,r)=B_{F}^{\alpha}(r,w)$. Since, we have
\begin{equation}\label{k981}
B_{F}^{\alpha}(r,w)=
\int_{S_{F}^{\alpha}(0)}^{S_{F}^{\alpha}(1)}(S_{F}^{\alpha}(x))^{r-1}
(1-S_{F}^{\alpha}(x))^{w-1}d_{F}^{\alpha}x,
\end{equation}
using the transformation $S_{F}^{\alpha}(x)=1-S_{F}^{\alpha}(y)$,
we conclude that
\begin{equation}\label{ds7}
   B_{F}^{\alpha}(r,w) =
\int_{S_{F}^{\alpha}(0)}^{S_{F}^{\alpha}(1)}(1-S_{F}^{\alpha}(y))^{r-1}
(S_{F}^{\alpha}(y))^{w-1}d_{F}^{\alpha}y=B_{F}^{\alpha}(w,r).
\end{equation}
2) Using the transformation
$S_{F}^{\alpha}(x)=\sin^{2}(S_{F}^{\alpha}(\theta))$, we get
following form for the fractal the Beta function
\begin{align}\label{nbvc}
&B_{F}^{\alpha}(r,w)=\int_{S_{F}^{\alpha}(0)}^{S_{F}^{\alpha}(\pi/2)}
\sin^{2}(S_{F}^{\alpha}(\theta))^{r-1}\cos^{2}(S_{F}^{\alpha}(\theta))^{w-1}
(2\sin(S_{F}^{\alpha}(\theta))\cos(S_{F}^{\alpha}(\theta)))d_{F}^{\alpha}x,\\&
    =2\int_{S_{F}^{\alpha}(0)}^{S_{F}^{\alpha}(\pi/2)}
    \sin^{2r-1}(S_{F}^{\alpha}(\theta))\cos^{2w-1}(S_{F}^{\alpha}(\theta))
    d_{F}^{\alpha}x.
\end{align}
3) The Beta fractal function is related to the fractal Gamma
 function as
\begin{equation}\label{vcxdzsax}
B_{F}^{\alpha}(r,w)=\frac{\Gamma_{F}^{\alpha}(r)\Gamma_{F}^{\alpha}(w)}
{\Gamma_{F}^{\alpha}(r+w)}.
\end{equation}
\textbf{Proof:} We have
\begin{equation}\label{mkiuj96}
    \Gamma_{F}^{\alpha}(r)\Gamma_{F}^{\alpha}(w)=4
    \int_{S_{F}^{\alpha}(0)}^{S_{F}^{\alpha}(\infty)}
    (S_{F}^{\alpha}(x))^{2r-1}(S_{F}^{\alpha}(y))^{2w-1}
    e^{S_{F}^{\alpha}(x)^2+S_{F}^{\alpha}(y)^2}d_{F}^{\alpha}x
    d_{F}^{\alpha}y.
\end{equation}
Transforming to polar coordinates
$S_{F}^{\alpha}(x)=S_{F}^{\alpha}(\rho)\cos(S_{F}^{\alpha}(\phi)),~~~
S_{F}^{\alpha}(x)=S_{F}^{\alpha}(\rho)\sin(S_{F}^{\alpha}(\phi))$
we obtain
\begin{align}\label{qwa}
    \Gamma_{F}^{\alpha}(r)\Gamma_{F}^{\alpha}(w)&=4\int_{S_{F}^{\alpha}(0)}^{S_{F}^{\alpha}(\infty)}
    (S_{F}^{\alpha}(\rho))^{2(r+w)-1}e^{-S_{F}^{\alpha}(\rho)^2}d_{F}^{\alpha}\rho
\int_{S_{F}^{\alpha}(0)}^{S_{F}^{\alpha}(\pi/2)}\cos^{2r-1}(S_{F}^{\alpha}(\phi))
\sin^{2w-1}(S_{F}^{\alpha}(\phi))d_{F}^{\alpha}\phi,\nonumber\\&
    =B_{F}^{\alpha}(r,w)\Gamma_{F}^{\alpha}(r+w).
\end{align}
Thus, the proof is completed.
\section{Non-local fractal derivative and integral \label{se-3}}

\begin{figure}[H]
\center
   \includegraphics[scale=0.6]{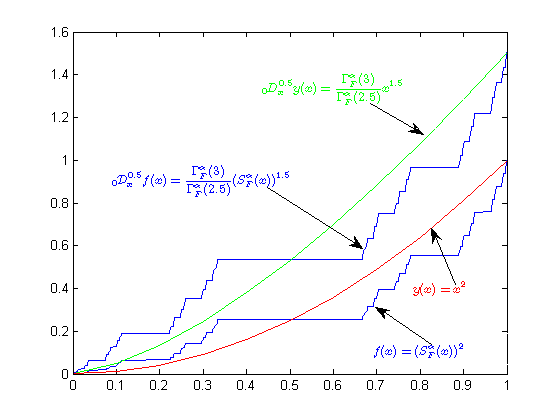}
 \caption{We plot  $y(x)=x^2$ and $f(x)=S_{F}^{\alpha}(x)^2$ and their non-local
 derivative   $_{0}D_{x}^{0.5}y(x)$ and $_{0}\mathcal{D}_{x}^{0.5}f(x)$, respectively.
 }\label{5414}
\end{figure}

\begin{figure}[H]
\center
   \includegraphics[scale=0.5]{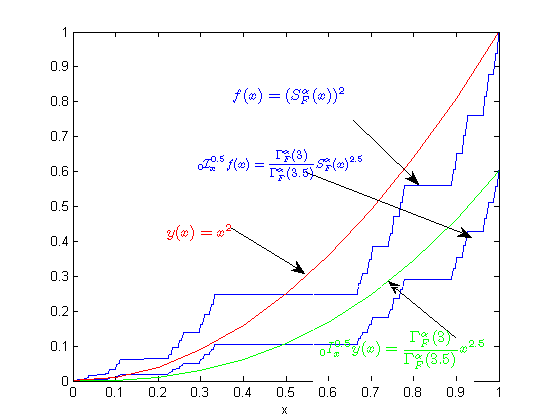}
 \caption{We show the graph of  $g(x)=x^2$ and $f(x)=S_{F}^{\alpha}(x)^2$ and
  their non-local integral    $_{0}I_{x}^{0.5}g(x)$ and $_{0}
  \mathcal{I}_{x}^{0.5}f(x)$, respectively.   }\label{544}
\end{figure}

In this section, we define the non-local derivative for the
functions
with fractal support.\\
 \textbf{Definition 3.} If $f(x)\in
C_{F}^{\alpha}[a,b]$ ($\alpha$-order differentiable function
 on $[a,b]$)  and $\beta>0$   then we have
\begin{equation}\label{xsaz125}
_{a}\mathcal{I}_{x}^{\beta}f(x):=\frac{1}{\Gamma_{F}^{\alpha}(\beta)}\int_{S_{F}^{\alpha}(a)}^
{S_{F}^{\alpha}(x)}\frac{f(t)}{(S_{F}^{\alpha}(x)-S_{F}^{\alpha}(t))^{\alpha-\beta}}
d_{F}^{\alpha}t,~~~ ~S_{F}^{\alpha}(x)>S_{F}^{\alpha}(a),
\end{equation}
where if $\beta=\alpha$ then we have  fractal integral whose order
 is equal the dimension of the fractal, and
\begin{equation}\label{frxsaz125}
_{x}\mathcal{I}_{b}^{\beta}f(x):=\frac{1}{\Gamma_{F}^{\alpha}(\beta)}
\int_{S_{F}^{\alpha}(x)}^
{S_{F}^{\alpha}(b)}\frac{f(t)}{(S_{F}^{\alpha}(x)-S_{F}^{\alpha}(t))^{\alpha-\beta}}
d_{F}^{\alpha}t,~~~ ~S_{F}^{\alpha}(x)<S_{F}^{\alpha}(b),
\end{equation}
are called the analogous left sided and the right sided
Riemann-Liouville fractal integral of
order $\beta$.\\
\textbf{Definition 4.} Let $n-\alpha\leq\beta < n$, then analogous
left and right Riemann-Liouville fractal derivative are defined as
follows:
\begin{align}\label{gfvgf}
 & _{a}\mathcal{D}_{x}^{\beta} f(x):=\frac{1}{\Gamma_{F}^{\alpha}(n-\beta)} (D_{F}^{\alpha})^n
\int_{S_{F}^{\alpha}(a)}^{S_{F}^{\alpha}(x)}\frac{f(t)}
{(S_{F}^{\alpha}(x)-S_{F}^{\alpha}(t))^{-n+\beta+\alpha}}d_{F}^{\alpha}t,\\&
 _{x}\mathcal{D}_{b}^{\beta} f(x):=\frac{1}{\Gamma_{F}^{\alpha}(n-\beta)} (-D_{F}^{\alpha})^n
\int_{S_{F}^{\alpha}(x)}^{S_{F}^{\alpha}(b)}\frac{f(t)}
{(S_{F}^{\alpha}(t)-S_{F}^{\alpha}(x))^{-n+\beta+\alpha}}d_{F}^{\alpha}t.
\end{align}
\textbf{Definition 5.} Let $f(x)\in C^{\alpha n} [a, b]$, then the
analogous left sided Caputo  fractal derivative is defined by
\begin{equation}\label{gfvsewgf}
  _{a}^{C}\mathcal{D}_{x}^{\beta} f(x):=\frac{1}{\Gamma_{F}^{\alpha}(n-\beta)}
\int_{S_{F}^{\alpha}(a)}^{S_{F}^{\alpha}(x)}(S_{F}^{\alpha}(x)-S_{F}
^{\alpha}(t))^{n-\beta-\alpha} (D_{F}^{\alpha})^nf(t)
d_{F}^{\alpha}t,~~~n=\max(0,-[-\beta]).
\end{equation}
Also, the analogous right sided Caputo fractal derivative has the
form
\begin{equation}
 _{x}^{C}\mathcal{D}_{b}^{\beta} f(x):=\frac{1}{\Gamma_{F}^{\alpha}(n-\beta)}
\int_{S_{F}^{\alpha}(x)}^{S_{F}^{\alpha}(b)}
(S_{F}^{\alpha}(t)-S_{F}^{\alpha}(x))^{n-\beta-\alpha}(-D_{F}^{\alpha})^n
f(t) d_{F}^{\alpha}t.
\end{equation}
Now, we give some important relations, namely
\begin{equation}\label{xsaza}
  _{a}\mathcal{I}_{x}^{\beta}(S_{F}^{\alpha}(x)-S_{F}^{\alpha}(a))^{\eta}=
  \frac{\Gamma_{F}^{\alpha}(\eta+1)}
{\Gamma_{F}^{\alpha}(\eta+\beta+1)}(S_{F}^{\alpha}(x)-S_{F}^{\alpha}(a))^{\eta+\beta},~~~\eta>-1.
\end{equation}
\textbf{Proof:} Using the Eq. (\ref{xsaz125}) we conclude
\begin{equation}\label{m258jyhnbgt}
_{a}\mathcal{I}_{x}^{\beta}(S_{F}^{\alpha}(x)-S_{F}^{\alpha}(a))^{\eta}=
\frac{1}{\Gamma_{F}^{\alpha}(\beta)}
\int_{S_{F}^{\alpha}(a)}^{S_{F}^{\alpha}(x)}
(S_{F}^{\alpha}(x)-S_{F}^{\alpha}(t))^{\beta-1}
(S_{F}^{\alpha}(t)-S_{F}^{\alpha}(a))^{\eta}d_{F}^{\alpha}t.
\end{equation}
Let us consider
\begin{equation}\label{bgthcxd963}
S_{F}^{\alpha}(\xi)=
\frac{S_{F}^{\alpha}(t)-S_{F}^{\alpha}(a)}{S_{F}^{\alpha}(x)-S_{F}^{\alpha}(a)},~~~~~~~~
d_{F}^{\alpha}t=(S_{F}^{\alpha}(x)-S_{F}^{\alpha}(a))d_{F}^{\alpha}\xi.~~~~~
\end{equation}
Therefore,  $S_{F}^{\alpha}(\xi):S_{F}^{\alpha}(0)\rightarrow
    S_{F}^{\alpha}(1)$ while $ S_{F}^{\alpha}(t):S_{F}^{\alpha}(a)\rightarrow
    S_{F}^{\alpha}(x)$. As a result we obtain
\begin{equation}\label{cdesx44}
    S_{F}^{\alpha}(x)-S_{F}^{\alpha}(t)
    =\frac{S_{F}^{\alpha}(1)-S_{F}^{\alpha}(\xi)}
    {S_{F}^{\alpha}(\xi)}(S_{F}^{\alpha}(t)-S_{F}^{\alpha}(0)).
\end{equation}
Substituting Eqs. ( \ref{bgthcxd963}) and  ( \ref{cdesx44}) in Eq.
(\ref{m258jyhnbgt}) we conclude that
\begin{align}\label{xswed9663}
   &_{a}\mathcal{I}_{x}^{\beta}(S_{F}^{\alpha}(x)-S_{F}^{\alpha}(a))^{\eta} =\frac{1}{\Gamma_{F}^{\alpha}(\beta)}\int_{S_{F}^{\alpha}(0)}^{S_{F}^{\alpha}(1)}
    (1-S_{F}^{\alpha}(\xi))^{\beta-1}S_{F}^{\alpha}(\xi)^{1-\beta}
    (S_{F}^{\alpha}(t)-S_{F}^{\alpha}(a))^{\beta+\eta-1}
    (S_{F}^{\alpha}(x)-S_{F}^{\alpha}(a))d_{F}^{\alpha}\xi,\nonumber\\&
    =\frac{1}{\Gamma_{F}^{\alpha}(\beta)}\int_{S_{F}^{\alpha}(0)}^{S_{F}^{\alpha}(1)}
    (1-S_{F}^{\alpha}(\xi))^{\beta-1}\left(\frac{S_{F}^{\alpha}(t)-S_{F}^{\alpha}(a)}
    {S_{F}^{\alpha}(x)-S_{F}^{\alpha}(a)}\right)^{1-\beta}
    (S_{F}^{\alpha}(t)-S_{F}^{\alpha}(a))^{\beta+\eta-1}
    (S_{F}^{\alpha}(x)-S_{F}^{\alpha}(a))d_{F}^{\alpha}\xi.
\end{align}
Then, we have
\begin{equation}\label{xswedftrg9663}
_{a}\mathcal{I}_{x}^{\beta}(S_{F}^{\alpha}(x)-S_{F}^{\alpha}(a))^{\eta}=\frac{(S_{F}^{\alpha}(x)-S_{F}^{\alpha}(a))^{\beta+\eta}}{\Gamma_{F}^{\alpha}(\beta)}
\int_{S_{F}^{\alpha}(0)}^{S_{F}^{\alpha}(1)}
(1-S_{F}^{\alpha}(\xi))^{\beta-1} (S_{F}^{\alpha}(\xi))^{\eta}
d_{F}^{\alpha}\xi.
\end{equation}
In view of Eq. (\ref{zaws}) we derive
\begin{equation}\label{xswedftrg9663}
_{a}\mathcal{I}_{x}^{\beta}(S_{F}^{\alpha}(x)-S_{F}^{\alpha}(a))^{\eta}=
\frac{(S_{F}^{\alpha}(x)-S_{F}^{\alpha}(a))^{\beta+\eta}}{\Gamma_{F}^{\alpha}(\beta)}
B_{F}^{\alpha}(\beta,\eta+1).
\end{equation}
Applying  Eq. (\ref{vcxdzsax}) we get
\begin{align}\label{xswedftrg9663}
    _{a}\mathcal{I}_{x}^{\beta}(S_{F}^{\alpha}(x)-S_{F}^{\alpha}(a))^{\eta}&=
    \frac{(S_{F}^{\alpha}(x)-S_{F}^{\alpha}(a))^{\beta+\eta}}{\Gamma_{F}^{\alpha}(\beta)}
\frac{\Gamma_{F}^{\alpha}(\beta)\Gamma_{F}^{\alpha}(\eta+1)}
{\Gamma_{F}^{\alpha}(\beta+\eta+1)},\nonumber\\&
    =\frac{\Gamma_{F}^{\alpha}(\eta+1)}{\Gamma_{F}^{\alpha}(\beta+\eta+1)}
(S_{F}^{\alpha}(x)-S_{F}^{\alpha}(a))^{\beta+\eta}.
\end{align}
Now, we consider following formula
\begin{equation}\label{xsamjuza}
  _{a}\mathcal{D}_{x}^{\beta}(S_{F}^{\alpha}(x)-S_{F}^{\alpha}(a))^{\eta}=
  \frac{\Gamma_{F}^{\alpha}(\eta+1)}
{\Gamma_{F}^{\alpha}(\eta+1-\beta)}(S_{F}^{\alpha}(x)-S_{F}^{\alpha}(a))^{\eta-\beta}.
\end{equation}
\textbf{Proof:} By rewriting  the Eq. (\ref{xsamjuza}) we get
\begin{equation}\label{u6xsamjuza}
  _{a}\mathcal{D}_{x}^{\beta}(S_{F}^{\alpha}(x)-S_{F}^{\alpha}(a))^{\eta}=
  (D_{F}^{\alpha})^{n}~_{a}\mathcal{I}_{x}^{n-\beta}
  (S_{F}^{\alpha}(x)-S_{F}^{\alpha}(a))^{\eta}.
\end{equation}
Utilizing the Eq. (\ref{xsaza}) we conclude
\begin{align}\label{m}
    _{a}\mathcal{D}_{x}^{\beta}(S_{F}^{\alpha}(x)-S_{F}^{\alpha}(a))^{\eta}&=\frac{\Gamma_{F}^{\alpha}(\eta+1)}
{\Gamma_{F}^{\alpha}(\eta+n-\beta+1)}
(D_{F}^{\alpha})^{n}(S_{F}^{\alpha}(x)-S_{F}^{\alpha}(a))^{\eta+n-\beta},\\&
    =\frac{\Gamma_{F}^{\alpha}(\eta+1)}
{\Gamma_{F}^{\alpha}(\eta-\beta+1)}
(D_{F}^{\alpha})^{n}(S_{F}^{\alpha}(x)-S_{F}^{\alpha}(a))^{\eta-\beta},~~~\eta>-1.
\end{align}
Now, we write some  important composition relations, namely
\begin{equation}\label{es789w8512}
 _{a}\mathcal{I}_{x}^{\beta}~_{a}\mathcal{D}_{x}^{\beta}f(x)=f(x)-\sum_{j=1}^{n}
\frac{(_{a}\mathcal{D}_{x}^{\beta-j}f(x))|_{(S_{F}^{\alpha}(a))}}
{\Gamma_{F}^{\alpha}(\beta+1-j)}(S_{F}^{\alpha}(x)-S_{F}^{\alpha}(a))^{\beta-j}.
\end{equation}
\textbf{Proof:} Using the definitions we get
\begin{align}\label{bz9}
_{a}\mathcal{I}_{x}^{\beta}~_{a}\mathcal{D}_{x}^{\beta}f(x)&=
\frac{1}{\Gamma_{F}^{\alpha}(\beta)}\int_{S_{F}^{\alpha}(a)}
^{S_{F}^{\alpha}(x)}(S_{F}^{\alpha}(x)-S_{F}^{\alpha}(t))^{\beta-1}
\mathcal{D}_{x}^{\beta}f(t)d_{F}^{\alpha}t\\&
=\frac{1}{\Gamma_{F}^{\alpha}(\beta+1)}D_{F}^{\alpha}\int_{S_{F}^{\alpha}(a)}
^{S_{F}^{\alpha}(x)}(S_{F}^{\alpha}(x)-S_{F}^{\alpha}(t))^{\beta}
(D_{F}^{\alpha})^{n}~
_{a}\mathcal{I}_{x}^{n-\beta}f(t)d_{F}^{\alpha}t
f(t)d_{F}^{\alpha}t.
\end{align}
Applying, n-times integration by part it leads to
\begin{align}\label{bgtrf}
    _{a}\mathcal{I}_{x}^{\beta}~_{a}\mathcal{D}_{x}^{\beta}f(x)&=D_{F}^{\alpha}~_{a}\mathcal{I}_{x}^{\beta+1-n}(
    _{a}\mathcal{I}_{x}^{n-\beta}f(x))-\sum_{k=1}^{n}
    \frac{(D_{F}^{\alpha})^{n-k}~
    _{a}\mathcal{D}_{x}^{\beta-n}f(x)|_{S_{F}^{\alpha}(a)}}{\Gamma_{F}^{\alpha}(\beta-k+1)}
    (S_{F}^{\alpha}(x)-S_{F}^{\alpha}(a))^{\beta-k},\nonumber\\&
    =f(x)-\sum_{k=1}^{n}~
    \frac{_{a}\mathcal{D}_{x}^{\beta-k}f(x)|_{S_{F}^{\alpha}(a)}}
    {\Gamma_{F}^{\alpha}(\beta-k+1)}(S_{F}^{\alpha}(x)-S_{F}^{\alpha}(a))^{\beta-k}.
\end{align}
The similar proof works for the following formulas
\begin{align}\label{eswnbv8512}
 &_{x}\mathcal{I}_{b}^{\beta}~_{x}\mathcal{D}_{b}^{\beta}f(x)=f(x)-\sum_{j=1}^{n}
\frac{(_{x}\mathcal{D}_{b}^{\beta-j}f(x))|_{(S_{F}^{\alpha}(b))}}
{\Gamma_{F}^{\alpha}(\beta+1-j)}(S_{F}^{\alpha}(b)-S_{F}^{\alpha}(x))^{\beta-j},\\&
 _{a}\mathcal{I}_{x}^{\beta}~_{a}^{C}\mathcal{D}_{x}^{\beta}f(x)=f(x)-\sum_{j=1}^{n}
\frac{((D_{F}^{\alpha})^j f(x))|_{(S_{F}^{\alpha}(a))}}
{\Gamma_{F}^{\alpha}(j+1)}(S_{F}^{\alpha}(x)-S_{F}^{\alpha}(a))^{j},\\&
 _{x}\mathcal{I}_{b}^{\beta}~_{x}^{C}\mathcal{D}_{b}^{\beta}f(x)=f(x)-\sum_{j=1}^{n}
\frac{((D_{F}^{\alpha})^j f(x))|_{(S_{F}^{\alpha}(b))}}
{\Gamma_{F}^{\alpha}(j+1)}(S_{F}^{\alpha}(b)-S_{F}^{\alpha}(x))^{j}.
\end{align}
In figures \ref{5414} and \ref{544} we  compared the non-local
standard derivative versus non-local fractal derivative and the
generalized fractal integral.
\section{Generalized  functions in the non-local calculus on the fractal subset of real-line \label{se-4} }

In this section, we suggest the mathematical tools for solving the
non-local fractal  differential equations.
\subsection{Gamma function on fractal subset of real line}

Now, we define the Gamma function for the fractal calculus that
will be used in non-local calculus on fractals.
\subsection{ Mittag-Leffler function on fractal subset of real-line}

It is well known that the exponential function has important role
in the theory of standard differential equation.~The generalized
exponential function is called the Mittag-Leffer function and
plays an important role
for fractional differential equations \cite{book-1}.\\
\textbf{Definition 6.} The generalized two parameter $\eta,\nu$
Mittag-Liffler function on fractal $F$ with $\alpha$-dimension is
defined as
\begin{equation}\label{d}
  E_{F,\eta,\nu}^{\alpha}(x)=\sum_{k=0}^{\infty}\frac{S_{F}^{\alpha}(x)^k}
  {\Gamma_{F}^{\alpha}(\eta k+\nu)},~~~\eta>0,~~\nu\in \Re.
\end{equation}
In the special case we have the following results, namely
\begin{align}
  E_{F,1,1}^{\alpha}(x)&= e^{S_{F}^{\alpha}(x)}, \\
  E_{F,1,2}^{\alpha}(x) &= \frac{e^{S_{F}^{\alpha}(x)-1}}{S_{F}^{\alpha}(x)}, \\
  E_{F,2,1}^{\alpha}(x) &= \cosh(S_{F}^{\alpha}(x)) ,\\
 E_{F,2,2}^{\alpha}(x) &= \frac{\sinh(S_{F}^{\alpha}(x))}{S_{F}^{\alpha}(x)}.
\end{align}
\subsection{Non-local Laplace transformation on fractal subset of real-line }

The Laplace transformation is a very useful tool for solving
standard linear differential equation with constant
coefficients.~The generalized Laplace transformation is
 applied to solve the fractional differential equations. Thus, in this section,
 we generalized the Laplace transformation for the function with fractal support
 which is utilized to solve the non-local differential equation on the fractal set \cite{book-1}. \\
\textbf{Definition 7.} Laplace transformation for the function
$f(x)$ is denoted by $F(s)$ and it is defined as
\begin{equation}\label{sewa}
  \mathcal{F}_{F}^{\alpha}(S_{F}^{\alpha}(s))=\mathcal{L}_{F}^{\alpha}[f(x)]=
  \int_{S_{F}^{\alpha}(0)}^{S_{F}^{\alpha}(\infty)}f(x)
e^{-S_{F}^{\alpha}(s)S_{F}^{\alpha}(x)}d_{F}^{\alpha}x.
\end{equation}
Now, we give the fractal Laplace transformation  of some
functions. If we define the fractal convolution of two function
$f(x)$ and $g(x)$ as follows:
\begin{equation}\label{zsewnja1485}
  f(x)*g(x)=\int_{S_{F}^{\alpha}(0)}^{S_{F}^{\alpha}(x)}
  f(S_{F}^{\alpha}(x)-S_{F}^{\alpha}(\tau))g(S_{F}^{\alpha}(\tau)) d_{F}^{\alpha}\tau,
\end{equation}
the fractal Laplace transformation of power function of
$S_{F}^{\alpha}(x)$ is
\begin{equation}\label{14nhytgb}
  \mathcal{L}_{F}^{\alpha}[S_{F}^{\alpha}(x)]
  \int_{S_{F}^{\alpha}(0)}^{S_{F}^{\alpha}(\infty)}S_{F}^{\alpha}(x)^{\beta}e^
  {-S_{F}^{\alpha}(s){S_{F}^{\alpha}(x)}}d_{F}^{\alpha}x=
  \frac{\Gamma_{F}^{\alpha}(1+\beta)}{s^{\beta+1}}.
\end{equation}
\textbf{Lemma:} The Laplace transformation of the non-local
fractal Riemann-Liouville integral is given by
\begin{equation}\label{xsewza}
  \mathcal{L}_{F}^{\alpha}[ _{0}\mathcal{I}_{x}^{\beta}f(x)]=
  \frac{\mathcal{F}_{F}^{\alpha}(S_{F}^{\alpha}(s))}
{S_{F}^{\alpha}(s)^{\beta}}.
\end{equation}
\textbf{Proof:} The Laplace transform of the fractal
Riemann-Liouville integral is
\begin{equation}\label{147}
   \mathcal{L}_{F}^{\alpha}[ _{0}\mathcal{I}_{x}^{\beta}f(x)] =
    \mathcal{L}_{F}^{\alpha}\left[ \frac{1}{\Gamma_{F}^{\alpha}(\beta)}\int_{S_{F}^{\alpha}(0)}^
{S_{F}^{\alpha}(x)}\frac{f(t)}{(S_{F}^{\alpha}(x)-S_{F}^{\alpha}(t))^{\alpha-\beta}}
d_{F}^{\alpha}t\right].
\end{equation}
Using the Eqs.(\ref{zsewnja1485}) and (\ref{14nhytgb}) we arrive
at
\begin{align}
   \mathcal{L}_{F}^{\alpha}[ _{0}\mathcal{I}_{x}^{\beta}f(x)]
   &=\frac{1}{\Gamma_{F}^{\alpha}(\beta)}\mathcal{F}_{F}^{\alpha}(S_{F}^{\alpha}(s))
  \mathcal{L}_{F}^{\alpha}[S_{F}^{\alpha}(x)^{\beta-1}],\nonumber
  \\&=
  \frac{1}{\Gamma_{F}^{\alpha}(\beta)}\mathcal{F}_{F}^{\alpha}(S_{F}^{\alpha}(s))
  \frac{\Gamma_{F}^{\alpha}(\beta)}{S_{F}^{\alpha}(s)^\beta},\nonumber\\&=
  \frac{\mathcal{F}_{F}^{\alpha}(S_{F}^{\alpha}(s))}
{S_{F}^{\alpha}(s)^{\beta}}.
\end{align}
The fractal Laplace transform of the non-local  fractal
Riemann-Liouville derivative of
 order $\beta\in [0,1)$ is given by
\begin{equation}\label{poi78}
  \mathcal{L}_{F}^{\alpha}\{  _{0}\mathcal{D}_{x}^{\beta}f(x),x,s \}=S_{F}^{\alpha}(s)^{\beta}
\mathcal{F}_{F}^{\alpha}(s)-\sum_{k=1}^{n}S_{F}^{\alpha}(s)^{n-k}~
_{0}\mathcal{D}_{x}^{\beta-n+k-1}f(x)|_{S_{F}^{\alpha}(0)},
\end{equation}
where $n=[\beta]+1$. The fractal Laplace transform of the
non-local  fractal Caputo derivative
 of order $\beta\in [0,1)$ is given by
\begin{equation}\label{poi78}
  \mathcal{L}_{F}^{\alpha}\{  _{0}^{C}\mathcal{D}_{x}^{\beta}f(x),x,s \}=(S_{F}^{\alpha}(s))^{\beta}
\mathcal{F}_{F}^{\alpha}(s)-\sum_{k=1}^{n}S_{F}^{\alpha}(s)^{\beta-k}~
_{0}\mathcal{D}_{x}^{k-1}f(x)|_{S_{F}^{\alpha}(0)}.
\end{equation}
where $n=\max(0,-[-\beta])$.
\section{Non-local fractal  differential equations \label{se-5}}

In this section, we solve  some illustrative examples.\\

\begin{figure}[H]
\center
   \includegraphics[scale=0.6]{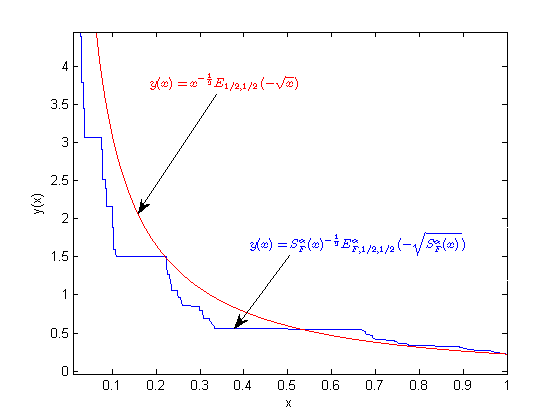}
 \caption{We plot the solution of Eq. (\ref{juyh}) on the real-line and Cantor set.
    }\label{8945lok44}
\end{figure}

\textbf{Example 1.} Consider the following linear fractal equation
\begin{equation}\label{xde}
  _{0}^{C}\mathcal{D}_{x}^{\frac{1}{2}}y(x)=2,
\end{equation}
with the initial condition
\begin{equation}\label{xz}
  D_{F}^{\alpha}y(x)|_{S_{F}^{\alpha}(0)=0}=1,
\end{equation}
where $\alpha= 0.6309$ is Cantor set dimension. By applying $
_{0}\mathcal{I}_{x}^{\frac{1}{2}}$ on the both sides of the  Eq.
(\ref{xz}) we obtain
\begin{equation}\label{bgtyhn}
  y(x)=\frac{1}{\Gamma_{F}^{\alpha}(1+\frac{1}{2})}S_{F}^{\alpha}(x)+
  \frac{2}{\Gamma_{F}^{\alpha}(1-\frac{1}{2})}S_{F}^{\alpha}(x)^{-\frac{1}{2}}.
\end{equation}
\begin{figure}[H]
\center
   \includegraphics[scale=0.6]{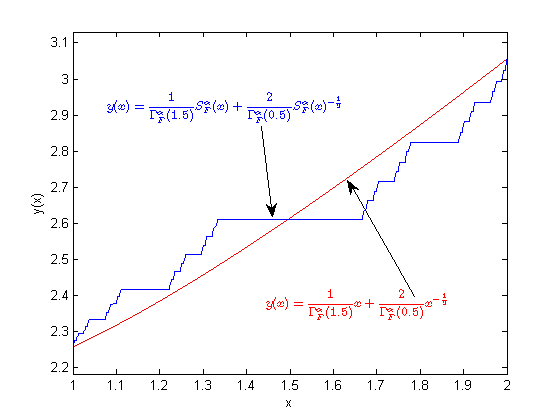}
 \caption{We present the solution of Eq. (\ref{xde}) on the real-line and Cantor set.
    }\label{8945lok4xxzsa4}
\end{figure}
\textbf{Example 2.} Consider a linear fractal differential
equation
\begin{equation}\label{swz5xc9}
   _{0}^{C}\mathcal{D}_{x}^{\frac{1}{2}}y(x)=1-S_{F}^{\alpha}(x),~~~~S_{F}^{\alpha}(x)\geq1,
\end{equation}
with initial condition as
\begin{equation}\label{n47712hytgb}
  D_{F}^{\alpha}y(x)|_{S_{F}^{\alpha}(1)}=0.
\end{equation}
By applying  $ _{0}\mathcal{I}_{x}^{\frac{1}{2}}$ on the both
sides of the Eq. (\ref{n47712hytgb})
 we arrive at
\begin{equation}\label{bwqalok28563}
  y(x)=-\frac{\Gamma_{F}^{\alpha}(2)}{\Gamma_{F}^{\alpha}(2+\frac{1}{2})}
  (S_{F}^{\alpha}(x)-1)^{1+\frac{1}{2}},~~~S_{F}^{\alpha}(x)\geq
  1.
\end{equation}
\begin{figure}[H]
\center
   \includegraphics[scale=0.6]{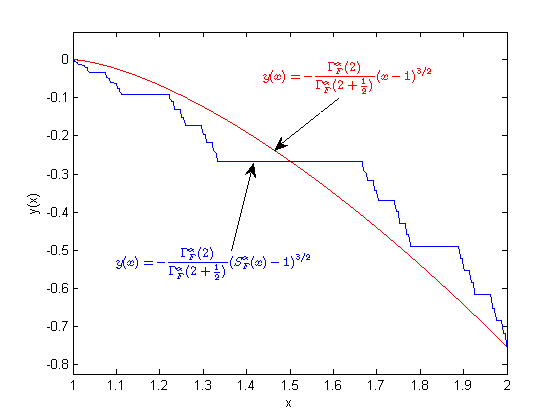}
 \caption{We give the graph of  the solution of Eq. (\ref{swz5xc9}) on the real-line and Cantor set.
    }\label{8945laqwqaqok44}
\end{figure}
In Figures, \ref{8945laqwqaqok44} and  \ref{8945lok4xxzsa4} we
plot the solutions of Eqs.(\ref{swz5xc9}) and (\ref{xde}),
respectively. \\
\textbf{Example 3.} Consider a linear differential equation
\begin{equation}\label{juyh}
  _{0}\mathcal{D}_{x}^{\frac{1}{2}}y(x)=~y(x),~~~~~~~~~~
\end{equation}
with the following initial condition, namely
\begin{equation}
_{0}\mathcal{D}_{x}^{-\frac{1}{2}}y(x)|_{S_{F}^{\alpha}(0)}=1.
\end{equation}
By inspection, the solution for the Eq. (\ref{juyh}) becomes
\begin{equation}\label{cxz}
  y(x)= S_{F}^{\alpha}(x)^{-\frac{1}{2}} E_{F,1/2,1/2}^{\alpha}(-\sqrt{S_{F}^{\alpha}(x)}).
\end{equation}
In  Figure \ref{8945lok44} we sketched the solution of  Eq.
(\ref{juyh}) on the Cantor set and real-line.
\textbf{Example 4.} We examine the following non-local
differential equation on a fractal subset of real-line, namely
with the following initial condition
\begin{align}\label{bds128}
  &_{0}\mathcal{D}_{F}^{\frac{4}{3}}y(x)- \lambda
  y(x)=(S_{F}^{\alpha}(x))^2,~\\&
_{0}\mathcal{D}_{F}^{\frac{1}{3}}y(x)|_{S_{F}^{\alpha}(0)}=1,~~~~~~~~~~~~~~~
_{0}\mathcal{D}_{F}^{\frac{-1}{6}}y(x)|_{S_{F}^{\alpha}(0)}=2.
\end{align}
For solving Eq. (\ref{bds128}) we apply the fractal Laplace
transformation on both side of it and we get
\begin{equation}\label{cfrdszw}
S_{F}^{\alpha}(s)^{\frac{4}{3}}\mathcal{F}_{F}^{\alpha}(s)-1-2(S_{F}^{\alpha}(s))^{\frac{1}{2}}
-\lambda
\mathcal{F}_{F}^{\alpha}(s)=\frac{2}{S_{F}^{\alpha}(s)^3}.
\end{equation}
After some calculations we obtain
\begin{equation}\label{n}
  \mathcal{F}_{F}^{\alpha}(s)=\frac{1}{S_{F}^{\alpha}(s)^{\frac{4}{3}}-\lambda}+
  \frac{2S_{F}^{\alpha}(s)^{\frac{1}{2}}}{S_{F}^{\alpha}(s)^{\frac{4}{3}}-\lambda}+
  \frac{2S_{F}^{\alpha}(s)^{-3}}{S_{F}^{\alpha}(s)^{\frac{4}{3}}-\lambda}.
\end{equation}
By computing the inverse fractal Laplace transform we conclude
\begin{align}\label{m}
  y(x)&=S_{F}^{\alpha}(x)^{\frac{4}{3}} E_{F,4/3,4/3}^{\alpha}
  (\lambda S_{F}^{\alpha}(x)^{\frac{4}{3}} )+2S_{F}^{\alpha}(x)^{\frac{-1}{6}}
  E_{F,4/3,5/6}^{\alpha}(\lambda S_{F}^{\alpha}(x)^{\frac{4}{3}}\nonumber\\&+2
   S_{F}^{\alpha}(x)^{\frac{10}{3}} E_{F,4/3,13/3}^{\alpha}(\lambda
   S_{F}^{\alpha}(x)^{\frac{4}{3}}.
\end{align}
\textbf{Remark:}~The figures \ref{8945lok4xxzsa4},
\ref{8945laqwqaqok44}, and \ref{8945lok44} show that the solution
of Eqs. (\ref{xde}),  (\ref{swz5xc9}) and (\ref{juyh}) leads to
the standard non-local fractional cases when  $\alpha=1$,
respectively.

\section{Conclusion \label{se-6}}

In this work, we defined new non-local derivatives on  fractal
sets.~These new type of non-local derivatives can  describe better
the dynamics of complex systems which possess memory effect on a
fractal set.~Four illustrative examples were solved in detail.
Finally, one can recover the standard non-local fractional cases
when put $\alpha=1$.

%%%%%%%%%%%%%%%%%%%%%%%%%%%%%%%%%%%%%%%%%%

%%%%%%%%%%%%%%%%%%%%%%%%%%%%%%%%%%%%%%%%%%

All authors common finished the manuscript. All authors have read
and approved the final manuscript.

%%%%%%%%%%%%%%%%%%%%%%%%%%%%%%%%%%%%%%%%%%

The authors declare no conflict of interest.

%=================================================================
% References: Variant A
%=================================================================
% Back Matter (References and Notes)
%----------------------------------------------------------
% Style and layout of the references

%=================================================================
% References:  Variant B
%=================================================================
% Use the following option to include external BibTeX files:
%\bibliography{lite}
%\bibliographystyle{mdpi}

%%%%%%%%%%%%%%%%%%%%%%%%%%%%%%%%%%%%%%%%%%

%\abbreviations{Abbreviations/Nomenclature}
%
%Main text.

%%%%%%%%%%%%%%%%%%%%%%%%%%%%%%%%%%%%%%%%%%

%\appendix
%\section{Appendix Title}
%
%Main text.

\end{document}